\documentclass[12pt, reqno]{amsart}
\usepackage{fullpage, amsthm, amssymb, stmaryrd, amsfonts, color, graphicx, hyperref}

\theoremstyle{plain} % default
\newtheorem*{thm*}{Theorem}
\newtheorem*{prob}{Problem}
\newtheorem{thm}{Theorem}[subsection]
\newtheorem{prop}[thm]{Proposition}
\newtheorem{lemma}[thm]{Lemma}
\newtheorem{cor}[thm]{Corollary}

\theoremstyle{definition}
\newtheorem{dfn}[thm]{Definition}

\newtheorem*{condition*}{Conditions}

\theoremstyle{remark}
\newtheorem{rem}[thm]{Remark}
\newtheorem*{rem*}{Remark}
\newtheorem*{question*}{Question}
\newtheorem*{example*}{Example}
\numberwithin{equation}{section}

\newcommand{\bb}{\mathbb}
\newcommand{\bbC}{\mathbb{C}}
\newcommand{\f}{\mathfrak}
\newcommand{\mc}{\mathcal}
\newcommand{\trdeg}{\mbox{trdeg }}

\newcommand{\la}{\lambda}
\newcommand{\q}{\mathbf{q}}
\renewcommand{\a}{\alpha}

\begin{document}
\title{Small semisimple subalgebras of \\ semisimple Lie algebras}
\author{Jeb F. Willenbring and Gregg Zuckerman} \date{\today}
\maketitle

\section{Introduction}

The main goal of this paper is to prove the following:
\begin{thm*}[see \cite{PZ}]  Let $\f k$ be an $\f{sl}_2$-subalgebra of a
semisimple Lie algebra $\f g$, none of whose simple factors is of
type $A1$. Then there exists a positive integer $b(\f k, \f g)$,
such that for every irreducible finite dimensional $\f g$-module
$V$, there exists an injection of $\f k$-modules $W\to V$, where
$W$ is an irreducible $\f k$-module of dimension less than $b(\f
k, \f g)$.
\end{thm*}

The goal of Section \ref{sec_Invariant_theory} is to provide a
proof of Theorem \ref{thm_invariants}. Section
\ref{sec_Representation_Theory} introduces the necessary facts
about Lie algebras and representation theory, with the goal being
the proof of Proposition \ref{prop_m_positive} (ultimately as an
application of Theorem \ref{thm_invariants}), and Proposition
\ref{prop_semigroup}.  In Section \ref{sec_Penkov_Zuckerman} we
prove the main theorem, using Propositions \ref{prop_semigroup}
and \ref{prop_m_positive}. In Section \ref{sec_G2}, we apply the
theorem to the special case where $\f g$ is the exceptional Lie
algebra $G_2$, and $\f k$ is a principal $\f{sl}_2$-subalgebra of
$\f g$.  We obtain a sharp estimate of $b(\f k, \f g)$ (in this
case).

\section{Invariant theory.}\label{sec_Invariant_theory}

The context for this section lies in the theory of algebraic group
actions on varieties. A good general reference for our terminology
and notation is \cite{AGIV} which contains translations of works,
\cite{PV} and \cite{Springer}.  For general notation and
terminology from commutative algebra and algebraic geometry see
\cite{Eisenbud} and \cite{Hartshorne}.  For the general theory of
linear algebraic groups, see \cite{Borel-LAG}.

All varieties are defined
over $\bbC$, although we employ some results that are valid in
greater generality.  Unless otherwise stated, all groups are assumed
to have the structure of a connected linear algebraic group. Of
particular interest is the situation where a group, $G$, acts on
an affine or quasi-affine variety. Given a (quasi-) affine
variety, $X$, we denote the ring of regular functions on $X$, by
$\bbC[X]$, and $\bbC[X]^G$ denotes the ring of $G$-invariant
functions.

We now turn to the following:
\begin{prob}
Let $X$ be an irreducible (quasi-) affine variety. Let $G$ be a
algebraic group acting regularly on $X$. When is $\bbC[X]^G \not
\cong \bbC$? That is to say, when do we have a non-trivial
invariant?
\end{prob}

\medskip The general problem may be too hard in this generality. We
begin by investigating a more restrictive situation which we
describe next.

\medskip For our purposes, a generic orbit, $\mc O \subset X$, is
defined to be an orbit of a point $x \in X$ with minimal isotropy
group (ie: $G_{x_0} = \{g \in G| g \cdot {x_0} = {x_0} \}$ for
$x_0 \in X$).
\begin{thm}\label{thm_invariants}
Assume that $X$ is an irreducible quasi-affine variety with a
regular action by a linear algebraic group $H$ such that:
\begin{enumerate}
    \item A generic $H$-orbit, $\mc O$, in $X$ has $\dim \mc O < \dim X$,
    \item $\bb C[X]$ is factorial, and
    \item $H$ has no rational characters\footnote{A rational character of $H$ is defined to be a regular function $\chi:H \rightarrow \bbC^{\times}$ such that $\chi(xy) = \chi(x)\chi(y)$ for all $x,y \in H$.}.
\end{enumerate}
Then, $\bbC[X]^H \not \cong \bb C$.
Furthermore, $\trdeg \bbC[X]^H = \mbox{codim } \mc O$.
\end{thm}

This theorem will be used to prove Proposition \ref{prop_invariants_general}.
In order to provide a proof of the above theorem we require some preparation.
Let $\bbC(X)$ denote the field of complex valued rational
functions on $X$. Our plan will be to first look at the ring of
rational invariants, $\bbC(X)^H$.  Given an integral domain, $R$,
let $Q R$ denote the quotient field. Clearly $Q \left( \bbC[X]^H
\right) \subseteq \bbC(X)^H$.  Under the assumptions of Theorem
\ref{thm_invariants} we have equality.  As is seen by:

\begin{thm}[\cite{AGIV} p. 165]\label{thm_AGIV_p165}
Suppose $k(X) = Q k[X]$.  If either,
\begin{itemize}
\item[(a)] the group $G^0$ is solvable\footnote{Notation: $G^0$ denotes the connected component of $G$.}, or
\item[(b)] the algebra $k[X]$ is factorial,
\end{itemize}
then any rational invariant of the action $G:X$ can be represented
as a quotient of two integral semi-invariants (of the same
weight). If, in addition, $G^0$ has no nontrivial characters
(which in case (a) means it is unipotent), then $k(X)^G = Q \left(
k[X]^G \right)$.
\end{thm}

As one might expect from the hypothesis (2) of Theorem
\ref{thm_invariants}, our applications of Proposition
\ref{thm_AGIV_p165} will involve condition (b).  We will then use:

\begin{prop}[\cite{AGIV} p. 166]\label{prop_AGIV_p166}
Suppose the variety $X$ is irreducible.  The algebra $k[X]^G$ separates orbits in general position
if and only if $k(X)^G = Q k[X]^G$, and in this case there exists a finite set of integral invariants that separates
orbits in general position and the transcendence degree of $k[X]^G$ is equal to the codimension
of an orbit in general position.
\end{prop}

\begin{proof}[Proof of Theorem \ref{thm_invariants}]
In our situation $X$ is a quasi-affine variety, thus $k(X) = Q
k[X]$.  By assumptions (2) and (3) of Theorem \ref{thm_invariants} and
Theorem \ref{thm_AGIV_p165} (using part (b)) we have that $k(X)^G
= Q \left( k[X]^G \right)$. $\dim \bbC[X]^H = \mbox{codim
} \mc O$ follows from the irreducibility of $X$ and Proposition
\ref{prop_AGIV_p166}. By assumption (1), $\mbox{codim } \mc O >
0$.  Thus, $\bbC[X]^H \not \cong \bbC$.
\end{proof}

\begin{question*}
For our purposes, we work within the context of the assumptions of
Theorem \ref{thm_invariants}, but to what extent may we relax the
assumptions to keep the conclusion of the theorem?
\end{question*}

\medskip \noindent Consider a triple $(G,S,H)$ such that the
following conditions ($\ast$) hold\footnote{Throughout this
article, we denote the Lie algebras of $G$ and $H$ by $\f g$ and
$\f k$, but we do not need this notation at present.}
\begin{condition*}[$\ast$]$\empty$
\begin{enumerate}
\item $G$ is a connected, simply-connected, semisimple linear algebraic group over $\bb C$.
\item $S$ and $H$ are connected algebraic subgroups of $G$ such that:
\begin{enumerate}
    \item $S \subseteq G$, is a connected algebraic subgroup with no non-trivial rational characters.
    \item $H \subseteq G$, is semisimple (and hence has no non-trivial rational characters).
\end{enumerate}
\end{enumerate}
\end{condition*}
\noindent
For our situation we will require the following two
results (used in Proposition \ref{prop_invariants_general}):
\begin{thm}[Voskresenski\u{\i} (see \cite{Popov-UFD}, \cite{Vosk-UFD},  \cite{Vosk})]\label{thm_simply-connected_factorial}
If $G$ is a connected, simply connected, semisimple linear
algebraic group then the ring of regular functions, $\bbC[G]$, is
factorial.
\end{thm}

\begin{thm}[see \cite{AGIV}, page 176]\label{thm_AGIV_p176}
If $\bbC[X]$ is factorial and the group $S$ is connected and has
no nontrivial characters, then $\bbC[X]^S$ is factorial.
\end{thm}

It is a theorem that $G/S$ has the structure of a quasi-affine
variety. Regarding the geometric structure of the quotient we
refer the reader to the survey in \cite{PV} Section 4.7. We
briefly summarize the main points: if $G$ is an algebraic group
with an algebraic subgroup $L \subseteq G$ then the quotient $G/L$
has the structure of a quasiprojective variety. If $G$ is
reductive, then $G/L$ is affine iff $L$ is reductive; $G/L$ is
projective iff $L$ is a parabolic subgroup (ie: $L$ contains a
Borel subgroup of $G$).  The condition of $G/L$ being quasi-affine
is more delicate, but includes, for example, the case when $L$ is
the first derived subgroup of a parabolic subgroup, which will be
of particular interest in Section \ref{subsec_parabolic}.

A consequence of Theorems \ref{thm_simply-connected_factorial} and \ref{thm_AGIV_p176} is that the algebra $\bbC[X]$ is factorial, where $X = G/S$.
$H$ then acts regularly on $G/S$ by left translation (ie:  $x \cdot gS = xgS$, for $g \in G$, $x \in H$).
And therefore, we are in a position to apply Theorem \ref{thm_invariants}, as in the next section.

\section{Representation Theory.}\label{sec_Representation_Theory}
This section begins by recalling some basic notation, terminology
and results of Lie theory.  We refer the reader to
\cite{Bourbaki-Lie123}, \cite{Bourbaki-Lie456}, \cite{FH},
\cite{GW}, and \cite{Knapp-Beyond} for this material.

\subsection{Notions from Lie theory.}\label{subsec_Lie_theory}$\empty$
Let $G$ denote a semisimple, connected, complex algebraic group.
We will assume that $G$ is simply connected. $T$ will denote a
maximal algebraic torus in $G$.  Let $r = \dim T$.  Let $B$ be a
Borel subgroup containing $T$. The unipotent radical, $U$, of $B$
is a maximal unipotent subgroup of $G$ such that $B = T \cdot U$.
Let $\f g$, $\f h$, $\f b$, $\f n^+$ denote the Lie algebras of
$G$, $T$, $B$, and $U$ respectively. Let $W := N_G(T)/T$ denote
the Weyl group corresponding to $G$ and $T$.

The Borel subalgebra, $\f b$, contains the Cartan subalgebra, $\f
h$ and is a semidirect sum\footnote{The signs $\subsetplus$ and
$\niplus$ stand for semidirect sum of Lie algebras.
If $\f g = \f g_1 \subsetplus \f g_2$
(or $\f g = \f g_2 \niplus \f g_1$) then $\f g_1$ is an ideal in $\f g$
and $\f g/\f g_1 \cong \f g_2$.}
$\f b = \f h \subsetplus \f n^+$.

The weights of $\f g$ are the linear functionals $\xi \in \f h^*$.
For $\a \in \f h^*$, set:
\[ \f g_\a=\{X\in\f g \; | \; [H,X]=\a(H) X \; \forall H \in \f h \}.\]
For $0 \neq \a \in \f h^*$, we say that $\a$ is a root if
$\f g_{\a} \neq (0)$.  For such $\a$, we have $\dim \f
g_\a = 1$.  Let $\Phi$ denote the set of roots. We then have
the decomposition:
\[\f g = \f h \oplus \sum_{\a \in \Phi} \f g_\a. \]
The choice of $B$ defines a decomposition $\Phi = \Phi^+ \cup
-\Phi^+$ so that $\f n^+ = \sum_{\a \in \Phi^+} \f g_\a$. We refer
to $\Phi^+$ (resp. $\Phi^- := - \Phi^+$) as the positive (resp.
negative) roots. Set: $\f n^- = \sum_{\a \in \Phi^+} \f g_{-\a}$.
Let $\overline B$ denote the (opposite) Borel subgroup of $G$ with
Lie algebra $\f h \oplus n^-$.  There is a unique choice of simple
roots $\Pi = \{ \a_1, \cdots, \a_r\}$ contained in $\Phi^+$, such
that each $\a \in \Phi^+$ can be expressed a non-negative integer
sum of simple roots.  $\Pi$ is a vector space basis for $\f h^*$.
Given $\xi, \eta \in \f h^*$ we write $\xi \preceq \eta$ if $\eta
- \xi$ is a non-negative integer combination of simple (equiv.
positive) roots. $\preceq$ is the dominance order on $\f h^*$.

For each positive root $\a$, we may choose a triple: $X_\a \in \f
g_\a$, $X_{-\a} \in \f g_{-\a}$ and $H_\a \in \f h$, such that
$H_\a = [X_\a, X_{-\a}]$ and $\a(H_\a) = 2$. Span $\{X_\a,
X_{-\a}, H_\a \}$ is then a three dimensional simple (TDS)
subalgebra of $\f g$, and is isomorphic to $\f{sl}_2$.

The adjoint representation, $ad: \f g \rightarrow \text{End}(\f
g)$ allows us to define the Killing form,
$(X,Y) = \text{Trace}(ad \, X \; ad \, Y)$ ($X,Y \in \f g$).
The semisimplicity of $\f g$
is equivalent to the non-degeneracy of the Killing form. By
restriction, the form defines a non-degenerate form on $\f h$,
also denoted $(,)$. Using this form we may define $\iota: \f h
\rightarrow \f h^*$ by $\iota(X)(-) = (X, -)$ ($X \in \f h$),
which allows us to identify $\f h$ with $\f h^*$.  Under this
identification, we have $\iota(H_\a) = \frac{2
\a}{(\a, \a)} =: \a^{\vee}$.

By definition, the Weyl group, $W$ acts on $T$.  By
differentiating this action we obtain an action on $\f h$, which
is invariant under $(,)$.  Via $\iota$, we obtain an action of $W$
on $\f h^*$.  In light of this, we view $W$ as a subgroup of the
orthogonal group on $\f h^*$.  $W$ preserves $\Phi$. For each
$\a \in \Phi$, set $s_\a(\xi) = \xi - (\xi, \a^\vee)
\a$ (for $\xi \in \f h^*$) to be the reflection through the
hyperplane defined by $\a^\vee$. We have $s_\a \in W$. For
$\a_i \in \Pi$, let $s_i := s_{\a_i}$, be the simple
reflection defined by $\a_i$. $W$ is generated by the simple
reflections. For $w \in W$, let $w = s_{i_1} s_{i_2} \cdot
s_{i_\ell}$ be a reduced expression (ie: an expression for $w$
with shortest length). The number $\ell$ is independent of the
choice of reduced expression.  We call $\ell =: \ell(w)$ the
length of $w$.  Note that $\ell(w) = |w(\Phi^+) \cap \Phi^-|$.
There is a unique longest element of $W$, denoted $w_0$ of length
$|\Phi^+|$.

The fundamental weights, $\{\omega_1, \cdots, \omega_r \}$ are
defined as $(\omega_i, \a_j^{\vee}) = \delta_{i,j}$.  We fix the
ordering of the fundamental weights to correspond with the usual
numbering of the nodes in the Dynkin diagram as in
\cite{Bourbaki-Lie456}.  Set $\rho = \frac{1}{2} \sum_{\a \in
\Phi^+} \a = \sum_{i=1}^r \omega_i$. A weight $\xi \in \f h^*$ is
said to be dominant if $(\la, \a) \geq 0$ for all $\a \in \Pi$.
The weight lattice $P(\f g) = \{\xi \in \f h^*| (\xi, \a^{\vee})
\in \bb Z \} = \sum_{i=1}^r \bb Z \omega_i$. We define the
dominant integral weights to be those $\xi \in P(\f g)$ such that
$(\xi, \a) \geq 0$ for all $\a \in \Pi$. The set of dominant
integral weights, $P_+(\f g)$, parameterizes the irreducible
finite dimensional representations of $\f g$ (or equivalently, of
$G$). We have \footnote{As usual, $\bb N = \{0,1,2, \cdots \}$
(the non-negative integers).} $P_+(\f g) = \sum_{i=1}^r \bb N
\omega_i$.

\subsection{Notions from representation theory.}\label{subsec_rep_theory}$\empty$
$\mc U(\f g)$ denotes the universal enveloping algebra of $\f g$.
The category of Lie algebra representations of $\f g$ is
equivalent to the category of $\mc U(\f g)$--modules.  A $\f
g$-representation (equiv. $\mc U(\f g)$--module), $M$, is said to
be a weight module if $M = \bigoplus M(\xi)$, where:
\[ M(\xi) = \left\{v \in M| Hv = \xi(H)v \; \forall H \in \f h \right\}. \]

Among weight modules are the modules admitting a highest weight
vector.  That is to say, a unique (up to scalar multiple) vector,
$v_0 \in M$ such that: \begin{enumerate} \item $\bbC v_0 = M^{\f
n^+} := \{v \in M| \f n^+ \, v = 0 \}$, \item $M(\la) = \bbC v_0$
for some $\la \in \f h^*$, and \item $\mc U(\f n^-)v_0 = M$.
\end{enumerate} Such a module is said to be a highest weight module (equiv
highest weight representation). $\la$ is the highest weight of
$M$. Given $\xi \in \f h^*$ with $M(\xi) \neq (0)$ we have $\xi
\preceq \la$.

For $\la \in \f h^*$, we let $\bbC_\la$ be the 1-dimensional
representation of $\f h$ defined by $\la$, then extended trivially
to define a representation of $\f b$ by requiring $\f n^+ \cdot
\bbC_\la = (0)$. Let $N(\la) := \mc U(\f g) \otimes_{\mc U(\f b)}
\bbC_\la$ denote the Verma module defined by $\la$. Let $L(\la)$
denote the irreducible quotient of $N(\la)$.
For $\la, \mu \in \f h^*$, $L(\la) \cong L(\mu)$ iff $\la = \mu$. $L(\la)$ (and
$N(\la)$) are highest weight representations.  Any irreducible
highest weight representation is equivalent to $L(\la)$ for a unique
$\la \in \f h^*$.  The theorem of the highest weight asserts that
$\dim L(\la) < \infty$ iff $\la \in P_+(\f g)$.

Each $\mu \in P(\f g)$, corresponds to a linear character of $T$,
denoted $e^\mu$.  For $\la \in P_+(\f g)$, the character of
$L(\la)$ defines a complex valued regular function on $T$.  This
character may be expressed as in the following:
\begin{thm}[Weyl]\label{thm_Weyl_character_formula}
\[
    ch \; L(\la) = \frac{\sum_{w \in W} (-1)^{\ell(w)} e^{w(\la + \rho)-\rho}}{\prod_{\a \in \Phi^+} (1 -
    e^{-\a})}.
\]
\end{thm}
It becomes necessary for us to refer to representations of both the group $G$ and the Lie algebra of $G$ (always denoted $\f g$).
By our assumptions on $G$ (in conditions ($\ast$)), we have that the every
finite dimensional complex representation of $\f g$ integrates to a regular representation of $G$.  The differential
of this group representation recovers the original representation of the Lie algebra.
We will implicitly use this correspondence.

For our purposes, of particular importance is the decomposition of the regular \footnote{We use the word ``regular'' in two senses, the other being in the context of algebraic geometry.}
representation of $G$.  That is,
\begin{thm}\label{thm_Peter_Weyl}
For $f \in \bbC[G]$, $(g,h) \in G \times G$ define: $(g,h) \cdot
f(x) = f(g^{-1}x h)$ ($x \in G$). Under this action we have the
classical Peter-Weyl decomposition:
\begin{equation}
    \bbC[G] \cong \bigoplus_{\la \in P_+(\f g)} L(\la)^* \; \widehat \otimes \; L(\la),
\end{equation}
as a representation of $G \times G$.  Here the superscript $*$ denotes the dual representation.
Note that $L(\la)^*$ is an irreducible highest weight representation of highest weight $-w_0(\la)$.
\end{thm}

We introduce the following notation:
\begin{dfn}
Let $\sigma:G \rightarrow GL(V)$ (resp. $\tau:H \rightarrow GL(W)$)
be a representation of a group $G$ (resp. $H$).
If $H \subseteq G$ we may regard $(\sigma,V)$ as a representation of $H$
by restriction.  We set:
\[ [V, W] := \dim \text{Hom}_H(V,W). \]
If either $V$ or $W$ is infinite dimensional $\text{Hom}_H(V,W)$ may be infinite
dimensional.  In this case $[V,W]$ should be regarded as an infinite cardinal.
(We will not encounter this situation in what is to follow.)
If $V$ is completely reducible as an $H$-representation,
and $W$ is irreducible, $[V,W]$ is the multiplicity of $W$ in $V$ (By Schur's lemma).

Note: We will use the same (analogous) notation of the category of Lie algebra representation.
\end{dfn}

\subsection{On the $\chi(T)$-gradation of $\bbC[U\backslash G]$}
We apply a philosophy taught to us by R. Howe.  As before, $U$ is a maximal
unipotent subgroup of $G$, $T$ a maximal torus (normalizing $U$).

As in Theorem \ref{thm_Peter_Weyl}, $U \times G \subseteq G \times
G$ acts on $\bbC[G]$.  We have $\bbC[G]^U \cong \bbC[U \backslash
G]$. As $T$ normalizes $U$ we have an action of $T$ on $\bbC[U
\backslash G]$ via $t \cdot f(x) = f(t^{-1}x)$ for $t \in T$
and $x \in U \backslash G$. We call this action the \emph{left}
action, since the multiplication is on the left. By Theorem
$\ref{thm_Peter_Weyl}$, we have:
\begin{equation}\bbC[U \backslash G] = \bigoplus_{\la \in P_+(\f g)} \left(L(\la)^*\right)^U \otimes L(\la). \end{equation}
A consequence of the theorem of the highest weight is that $\dim \left(L(\la)^*\right)^U = 1$.
We let, $\chi(T) \cong \bb Z^r$, denote the character group of $T$.  Each $\la \in P_+(\f g)$
defines a character, $e^\la$, of $T$.
Set: $\bbC[U \backslash G]_\la := \left\{f \in \bbC[U \backslash G] \; |  \; f(t^{-1}x) = e^\la(t)f(x) \; \; \forall x \in U\backslash G, t \in T \right\}$.
$G$ then acts (by right multiplication) on $\bbC[U \backslash G]$, and under this action we have:
$\bbC[U \backslash G]_{\la} \cong L(\la)$.
We then obtain a $\chi(T)$-gradation of the algebra.  That is to say,
$\bbC[U \backslash G]_{\xi} \cdot \bbC[U \backslash G]_{\eta} \subseteq \bbC[U \backslash G]_{\xi + \eta}$.
We exploit this phenomenon to obtain:

\begin{prop}\label{prop_semigroup}
    Let $W$ be an irreducible finite dimensional representation of
    a reductive subgroup, $H$, of $G$.
    Given $\la, \mu \in P_+(\f g)$,
    If $L(\mu)^H \neq (0)$ and $[W, L(\la)] \neq 0$ then
    $[W, L(\la + \mu)] \neq 0$.
\end{prop}

\begin{rem}
    In Proposition \ref{prop_semigroup}, the $G$-representations $L(\la)$ and $L(\la + \mu)$
    are regarded as $H$-representations by restriction.
\end{rem}

\begin{proof}[Proof of Proposition \ref{prop_semigroup}]
Let $f \in \bbC[U \backslash G]_\mu^H \cong L(\mu)^H$ and $\tilde W \subseteq \bbC[U \backslash G]_\la \cong L(\la)$ such that $W \cong \tilde W$ as a representation of $H$.
Then $f \cdot W \subseteq \bbC[U \backslash G]_{\la + \mu}$.  Under the (right) action of $H$ we have, $f \cdot \tilde W \cong W$.
Therefore, $[W, \bbC[U\backslash G]_{\la + \mu}] \neq 0$.
\end{proof}

\subsection{The maximal parabolic subgroups of $G$.}\label{subsec_parabolic} $\empty$
A connected algebraic subgroup, $P$, of $G$ containing a Borel
subgroup is said to be parabolic.  There exists an inclusion
preserving one-to-one correspondence between parabolic subgroups
and subsets of $\Pi$.  We will recall the basic set-up.

Let $\f p = Lie(P)$ denote the Lie algebra of a parabolic subgroup
$P$. Then $\f p = \f h \oplus \sum_{\a \in \Gamma} \f g_\a$, where
$\Gamma := \Phi^+ \cup \{\a \in \Phi \; | \; \a \in
\text{Span}(\Pi^\prime) \}$ for a unique $\Pi^\prime \subseteq
\Pi$. Set:
\[\f l = \f h \oplus \sum_{\a \in \Gamma \cap -\Gamma} \f g_\a,
\;\;\; \f u^+ = \sum_{\a \in \Gamma, \, \a \not \in -\Gamma} \f
g_\a, \text{ and } \f u^- = \sum_{\a \in \Gamma, \, \a \not \in
-\Gamma} \f g_{-\a}.\] Then we have, $\f p = \f l \oplus \f u^+$
and $\f g = \f u^- \oplus \f p $. The subalgebra $\f l$ is the
Levi factor of $\f p$, while $\f u^+$ is the nilpotent radical of
$\f p$.  $\f l$ is reductive and hence $\f l = \f l_{ss} \oplus \f
z(\f l)$, where $\f z(\f l)$ and $\f l_{ss}$ denote the center and
semisimple part of $\f l$ respectively.

The following result is a slight modification of Exercise 12.2.4
in \cite{GW} (p. 532).

\begin{prop}\label{prop_GW_exercise}
For $0 \neq \la \in P_+(\f g)$, let $v_\la$ be a highest weight
vector in
$L(\la)$. \\
Let $X = G \cdot v_\la \subseteq L(\la)$ denote the orbit of
$v_\la$ and let $G_{v_\la} = \{g \in G| g \cdot v_{\la} = v_{\la}
\}$ denote the corresponding isotropy group.  Then, $X$ is a
quasi-affine variety stable under the action of $\bbC^\times$ on
$L(\la)$ defined by scalar multiplication.  This
$\bbC^\times$-action defines a graded algebra structure on $\bb
C[X] = \bigoplus_{d=0}^\infty \bbC[X]^{(d)}$.

The action of $G$ on $\bbC[X]$, defined by $g \cdot f(x) =
f(g^{-1} \cdot x)$ (for $g \in G$ and $x \in X$) commutes with the
$\bbC^\times$-action and therefore each graded component of
$\bbC[X]$ is a representation of $G$.  Furthermore, we have:
\[ \bbC[X]^{(n)} \cong L(-n w_0(\la))\]
for all $n \in \bb N$.
\end{prop}
\begin{proof} Set $V = L(\la)$, and let $v_\la$ be a highest weight vector in
$V$.  Let $\bb P V = \{[v]|v\in V \}$ denote
the complex projective space on $V$. $G$ then acts on $\bb P V$ by
$g \cdot [v] = [g \cdot v]$ for $g \in G$ and $v \in V$.  The
isotropy group of $[v_\la]$ contains the Borel subgroup, $B$ and
therefore is a parabolic subgroup, $P \subseteq G$.  $G/P$ is
projective so, $G \cdot [v_\la] \subset \bb PV$ is closed.
This means that the affine cone, $\bb A := \bigcup_{g \in G} [g \cdot v_\la]
\subset V$ is closed.  Let $a \in \bb A$. Then, $a = z (g \cdot
v_\la)$ for some $z \in \bbC$ and $g\in G$. The action of $G$ on
$V$ is linear so, $a = g \cdot zv_\la$.  By the assumption that
$\la \neq 0$, $T$ acts on $\bbC v_\la$ by a non-trivial linear
character. All non-trivial linear characters of tori are
surjective (this is because the image of a (connected) algebraic group
homomorphism is closed and connected).
Therefore, if $z \neq 0$, then $z v_\la = t \cdot
v_\la$ for some $t \in T$.  This fact implies\footnote{If $a \neq
0$, $a = g \cdot (t \cdot v_\la) = (gt) \cdot v_\la \in X.$} that
either $a \in X$ or $a = 0$. And so, $\overline X = X \cup \{0\} =
\bb A$. $X$ is therefore quasi-affine since $X = \bb A  - \{0\}$
($0 \not\in X$ since $\la \neq 0$).  We have also shown that $X$
is stable under scalar multiplication by a non-zero complex
number.

Let $v^*_\la \in L(\la)^*$ be a highest weight vector. Upon
restriction $v^*_\la$ defines a regular function on $X$ in
$\bbC[X]^{(1)}$, which is a highest weight vector for the left
action of $G$ on $\bbC[X]$.  Then, $(v^*_\la)^n \in
\bbC[X]^{(n)}$. $(v^*_\la)^n$ is a highest weight vector, hence,
$[\bbC[X]^{(n)}, L(n\la)^*] \neq 0$. That is to say, we have an
injective $G$-equivariant map, $\psi: L(n\la)^* \hookrightarrow
\bbC[X]^{(n)}$. It remains to show that $\psi$ is an isomorphism.

Since $X$ is quasi-affine, $\bbC[X] \cong \bbC[G]^{G_{v_\la}}$. By
restriction of the regular representation, $G \times G_{v_\la}$
acts on $\bbC[G]$.  And so by, Theorem \ref{thm_Peter_Weyl}:
\[
    \bbC[X] \cong \bigoplus_{\xi \in P_+(\f g)} L(\xi)^* \otimes L(\xi)^{G_{v_\la}}
\]
Set $\mc L_\xi := L(\xi)^{G_{v_\la}}$.  $U \subseteq G_{v_\la}$
since $v_\la$ is a highest weight vector. Therefore, $\mc L_\xi
\subseteq L(\xi)^U$ and $\dim \mc L_\xi \leq 1$.  And so, $\bb
C[X]$ is multiplicity free as a representation of $G$ (under the
action of left multiplication).\footnote{Alternatively, $\overline
B$ has a dense orbit in $X$ (since $\overline B U$ is dense in $G$
(ie: the big cell)).  A dense orbit under a Borel subgroup is
equivalent to having a multiplicity free coordinate ring.} We will
show that the only possible $\xi$ for which $\dim \mc L_\xi > 0$
are the non-negative integer multiples of $\la$. By the fact that
$\bbC[X]$ is multiplicity free we will see that $\psi$ must be an
isomorphism.

If $\mc L_\xi \neq (0)$ then $\mc L_\xi = L(\xi)^U$ since they are
both 1-dimensional. Assume that $\mc L_\xi \neq (0)$. Choose, $0
\neq v_\xi \in \mc L_\xi$. Note that $v_\xi$ is a highest weight
vector. $T$ acts on $\bbC v_\xi$ by $t \cdot v_\xi = e^\xi(t)
v_\xi$ for all $t \in T$. Set $T_\xi := T \cap G_{v_\la}$ and $\f
h_\la := Lie(T_\la) = \{H \in \f h \;|\; \la(H) = 0 \}$.

$H \in \f h_\la$ implies both $H \cdot v_\xi = 0$ and $H \cdot
v_\xi = \xi(H) v_\xi$.  Hence we have $\xi(H) = 0$ when $H \in \f
h_\la$.  This statement is equivalent to $\la$ and $\xi$ being
linearly dependent.  Furthermore, we have $n_1 \xi = n_2 \la$ for
$n_1, n_2 \in \bb N$ since $\xi$ and $\la$ are both dominant
integral weights.

If $\dim \mc L_\xi=1$ then $[\bbC[X]^{(n)},L(\xi)^*]\neq 0$ for
some $n \in \bb N$. This forces $[\bbC[X]^{(n n_1)}, L(n_1 \xi)^*]
\neq 0$ and therefore, $[\bbC[X]^{(n n_1)}, L(n_2 \la)^*] \neq 0$.
As before, $[\bbC[X]^{(n_2)}, L(n_2 \la)^*] \neq 0$.  Using the
fact that $\bbC[X]$ is multiplicity free, we have $n n_1 = n_2$.
And so, $\xi = n \la$.
\end{proof}

\begin{rem*}
The closure of the variety $X$ in Proposition
\ref{prop_GW_exercise} is called the \emph{highest weight variety}
in \cite{VP1972} (see \cite{GW}).
\end{rem*}

If $\Pi^\prime = \Pi - \{\a\}$ for some simple root $\a$, then the
corresponding parabolic subgroup is maximal (among proper parabolic
subgroups). Consequently, the maximal parabolic subgroups of $G$
may be parameterized by the nodes of the Dynkin diagram,
equivalently, by fundamental weights of $G$. Set: $\la =
-w_0(\omega_k)$. Let $v_k \in L(\la)^U$ be a highest weight
vector.  Define:
\[
X^{(k)}_{\f g} := G \cdot v_k \subseteq L(\la) \hspace{1cm}(1\leq k \leq r),
\] the orbit of $v_k$ under the action of $G$.  (When there is no chance of
confusion, we write $X^{(k)}$ for $X^{(k)}_{\f g}$.)

We have seen that this orbit has the structure of a quasi-affine
variety.  It is easy to see that the isotropy group, $S^{(k)}$ of
$v_k$ is of the form $S^{(k)} = [P_k,P_k]$, where $P_k$ denotes a
maximal parabolic subgroup of $G$.\footnote{This fact in not true
for the orbit of an arbitrary dominant integral weight.  In
general, $[\f p, \f p] \neq Lie(G_{v_\la})$, for any parabolic
subalgebra $\f p$.} Let $\f p := Lie(P_k)$ and $\f p = \f l \oplus
\f u^+$ denote the Levi decomposition.  We have $\f p = \f l +
n^+$ and $n^+ = [\f b, \f b] \subset [\f p, \f p]$ and $\f l_{ss}
= [\f l, \f l] \subseteq [\f p, \f p]$. Therefore, $[\f p, \f p] =
\f l_{ss} + n^+ = \f l_{ss} \oplus u^+$.  (Note that all that is
lost is $\f z(\f l)$.)  We say that $P_k$ is the parabolic
subgroup corresponding to the fundamental weight $\omega_k$.

We have $\dim X^{(k)} = \dim \f g / [\f p, \f p]$ and $2 \dim \f
g/[\f p, \f p] = \dim (\f g/\f l_{ss}) + 1$ (Note that $\dim \f
z(\f l) = 1$ since $P_k$ is maximal.) In light of these facts, we
see that the dimension of $X^{(k)}$ may be read off of the Dynkin
diagram.  The dimension is important for the proof of Corollary
\ref{cor_SL2_invariants_noA2}.  In Section \ref{sec_tables}, we
explicitly compute $\dim \f g/ \f l_{ss}$ for the exceptional Lie
algebras and low rank classical Lie algebras.

\begin{cor}\label{cor_G_decomposition}
For $X = X^{(k)}$ where $1 \leq k \leq r$ we have:
\begin{equation}\label{eqn_G_decomposition}
    \bbC[X] \cong \bigoplus_{n=0}^\infty L(n \omega_k)
\end{equation}
as a representation of $G$.
\end{cor}
\begin{proof}
The result is immediate from Proposition \ref{prop_GW_exercise}.
\end{proof}

\subsection{A consequence of Theorem \ref{thm_invariants}.}\label{subsec_consequence}$\empty$\\
The goal of Section \ref{sec_Invariant_theory} was to prove Theorem \ref{thm_invariants}.
We now apply this theorem to obtain the following:
\begin{prop}\label{prop_invariants_general}
Assume $G$ and $H$ satisfy conditions ($\ast$), and
we take $S = S^{(k)}$, for some $1 \leq k \leq r$. Set: $X := X^{(k)}$.
\[ \dim H < \dim X \implies \bbC[X]^H \not \cong \bbC. \]
\end{prop}
\begin{proof}
    If $\mc O$ is a generic $H$-orbit in $X$ then $\dim \mc O \leq \dim H < \dim X$.
    By Theorems \ref{thm_simply-connected_factorial} and \ref{thm_AGIV_p176},
    the algebra $\bbC[X]$ is factorial, because $\bbC[X] = \bbC[G]^S$.
    The result follows from Theorem \ref{thm_invariants}.
\end{proof}
We now provide a representation theoretic interpretation of this proposition as it relates to $\f{sl}_2$..
\subsubsection{\bf The $\f{sl}_2$-case.}\label{subsec_consequence_sl2}
Consider a triple $(G,S,H)$ such that the following
conditions ($\ast \ast$) hold:
\begin{condition*}[$\ast \ast$]$\empty$
\begin{enumerate}
\item $(G,S,H)$ satisfy conditions ($\ast$), and:
\item $\f g$ has no simple factor of Lie type A1.
\item $S = S^{(k)}$, where $1 \leq k \leq r$.  Set: $X := X^{(k)} \;(:=G/S)$. \item $\f{sl}_2
\cong Lie(H) =: \f k \subseteq \f g$
\end{enumerate}
\end{condition*}
\begin{cor}\label{cor_SL2_invariants_noA2}
    Assume conditions ($\ast \ast$), and that $\f g$ has no simple factor of Lie type $A2$.
    Then, $\bbC[X]^H \not \cong \bbC$.
\end{cor}
\begin{proof}
    The statement can be reduced to the case where $\f g$ is simple.
    For $G$ simple and not of type A1 or A2,
    we appeal to the classification of maximal parabolic subgroups
    (see the tables in Section \ref{sec_tables})
    to deduce that $\dim X > 3$.
    as $H \cong SL_2(\bbC)$ (locally).  Hence, $\dim H < \dim X$.
    We are then within the hypothesis of Proposition \ref{prop_invariants_general}.
\end{proof}
We next address the case when $\f g$ does have a simple factor of
Lie type $A2$.  For this material we need to analyze the set of
$\f{sl}_2$-subalgebras of $\f g = \f{sl}_3$, up to a Lie algebra
automorphism.  For general results on the subalgebras of $\f g$ we
refer the reader to \cite{Dynkin-maximal} and
\cite{Dynkin-semisimple}.

In the case of $\f g = \f{sl}_3$ there are two such $\f{sl}_2$
subalgebras.  One being the root $\f sl_2$-subalgebra
corresponding to any one of the three positive roots of $\f g$.
The other is the famous principal $\f sl_2$-subalgebra.

A principal $\f {sl}_2$-subalgebra (\cite{Dynkin-semisimple},
\cite{Kostant}) of $\f g$ is a subalgebra $\f k \subseteq \f g$ such that $\f k \cong {sl}_2$ and contains a
regular nilpotent element.  These subalgebras are conjugate, so we
sometimes speak of ``the'' principal $\f {sl}_2$-subalgebra.
There is a beautiful connection between the principal $\f{sl}_2$-subalgebra and the cohomology of $G$, (see \cite{Kostant}).
There is a nice discussion of this theory in \cite{CM}.

\begin{lemma}\label{lemma_A2_principal_sl2} Assume conditions ($\ast \ast$).
    Let $G$ and $H$ be such that $\f g \cong \f{sl}_3$, and
    $Lie(H)$ is a principal $\f{sl}_2$-subalgebra of $\f g$.
    Then, $\dim \bbC[X]^H > 0$, for $X=X^{(1)}$ or $X = X^{(2)}$.
\end{lemma}
\begin{proof}
A principal $\f {sl}_2$-subgroup in $\f {sl}_3$ is embedded as a symmetric subalgebra.
More precisely, let $H := SO_3(\bbC) \subseteq G$, then $Lie(H)$ is a principal $\f{sl}_2$-subalgebra of $\f {sl}_3$.

In general, if $G = SL_n(\bbC)$ and $K = SO_n(\bbC)$ with $H$
embedded in $G$ in the standard way, then the pair $(G,H)$ is
symmetric (ie: $H$ is the fixed point set of a regular involution
on $G$).  As before, let $r$ denote the rank of $G$. In order the
prove the lemma (for $r=2$), it suffices to observe that by the
Cartan-Helgason theorem (see \cite{GW}, Chapter 11) we have:

\[\dim L(\la)^K =
\left\{%
\begin{array}{ll}
    1, & \hbox{$\la \in \sum_{i=1}^r 2\bb N \omega_i$;} \\
    0, & \hbox{Otherwise.} \\
\end{array}%
\right.
\]
\end{proof}
\begin{lemma}\label{lemma_A2_root_sl2} Assume conditions ($\ast \ast$).
    Let $G$ and $H$ be such that $\f g \cong \f{sl}_3$, and
    $Lie(H)$ is a root $\f{sl}_2$-subalgebra of $\f g$.
    Then, $\dim \bbC[X]^H > 0$, for $X=X^{(1)}$ or $X = X^{(2)}$.
\end{lemma}
\begin{proof}
It is the case that $L(\omega_1)$ and $L(\omega_2)$ both have
$H$-invariants, as they are equivalent to the standard
representation of $G$ and its dual respectively.
\end{proof}
Summarizing we obtain:
\begin{prop}\label{prop_SL2_invariants}
    Under assumptions ($\ast \ast$), $\dim \bbC[X]^H > 0$.
\end{prop}
\begin{proof}
The statement reduces to the case where $\f g$ is simple, because a maximal parabolic subalgebra
of $\f g$ must contain all but one simple factor of $\f g$.  Therefore, assume that $\f g$ is simple
without loss of generality.
For $\f g \not \cong \f{sl}_3$, apply Corollary \ref{cor_SL2_invariants_noA2}.
If $\f g \cong \f{sl}_3$, then $Lie(H)$ is embedded as either a root $\f{sl}_2$-subalgebra,
or as a principal $\f{sl}_2$-subalgebra.  Apply Lemmas \ref{lemma_A2_root_sl2} and
\ref{lemma_A2_principal_sl2} to the respective cases.
\end{proof}

The following will be of fundamental importance in Section \ref{sec_Penkov_Zuckerman}.
\begin{dfn}\label{dfn_m}  For $G$ and $H$ as in conditions ($\ast$),
we consider the following set of positive integers \footnote{As
always, $\bb Z^+ = \{1, 2, 3, \cdots \}$ (the positive
integers).}:
\[ M(G,H, j) := \left \{ n\in \bb Z^+ | \dim \left[ L( n \omega_j) \right]^H \neq (0) \right \} \]
where $j$ is a positive integer with $1 \leq j \leq r$.  Set:
\[
m(G,H,j) :=
\left\{%
\begin{array}{cl}
    \min M(G, H, j), & \hbox{if $M(G,H,j)\neq \emptyset$;} \\
    0              , & \hbox{if $M(G,H,j)   = \emptyset$.} \\
\end{array}%
\right.
\]
We will also write $m(\f g, \f k, j)$ (resp. $M(\f g, \f k,
j$)) where (as before) $\f k = Lie(H)$ and $\f g = Lie(G)$.
\end{dfn}

\begin{prop}\label{prop_m_positive}
For $G$ and $H$ as in conditions ($\ast \ast$),
\[m(G,H, k)>0 \text{ for all, } 1 \leq k \leq r.\]
\end{prop}
\begin{proof}
    Apply Corollary \ref{cor_G_decomposition} and Proposition \ref{prop_SL2_invariants}.
\end{proof}
Proposition \ref{prop_invariants_general} applies to a much more general situation than $\f k \cong \f{sl}_2$.
\subsubsection{\bf The semisimple case.}
As it turns out, what we have done for the $\f{sl}_2$-subalgebras can be
done for any semisimple subalgebra.

\begin{prop}\label{prop_m_positive_general}
If $G$ and $H$ are as in condition ($\ast$), then for each $1 \leq k \leq r$,
\[\dim H < \dim X^{(k)} \implies m(G,H,k)>0.\]
\end{prop}
\begin{proof}
Apply Proposition \ref{prop_invariants_general} and Corollary \ref{cor_G_decomposition}.
\end{proof}

In order to effectively apply Proposition \ref{prop_m_positive_general} we will want
to guarantee that $\dim H < \dim X^{(k)}$ for all $1 \leq k \leq r$.  This will happen
if all simple factors of $\f g$ have sufficiently high rank.  This idea
motivates the following definition.  Consider a group $G$ (as in condition ($\ast$)) and define:
\[ e(\f g) := \min_{ 1 \leq k \leq r} \dim X^{(k)}_{\f g}  \]
For a semisimple complex Lie algebra $\f k$, define\footnote{Explanation for notation: $E$ and $e$ are chosen with the word ``exclusion'' in mind.}:
\[
    E(\f k) = \left\{\f s \; \left| \;
\begin{array}{l}
(1) \; \; \f s \text{ is a simple complex Lie algebra} \\
(2) \; \; \dim \f k \geq e(\f s)
\end{array}\right.
\right \}.
\]
\begin{cor}\label{cor_m_positive_general}
Assume $G$ and $H$ as in conditions ($\ast$).
If $\f g$ has no simple factor that is in the set $E(\f k)$
then $m(G, H, k) > 0$ for all $k$ with $1 \leq k \leq r$.
\end{cor}
\begin{proof}
    Immediate from Proposition \ref{prop_m_positive_general} and the definition of $E(\f k)$.
\end{proof}

\begin{example*}
\[    E(\f{sl}_3) = \{A_1, A_2, A_3 (=D_3), A_4, A_5, A_6, A_7, B_2 (=C_2), B_3, B_4, C_3, C_4, D_4, G_2 \} \]
By the tables in Section \ref{sec_tables}, we can determine:
\[
\begin{array}{c|cccccccccccccccccc}
\f s   & A_1 & A_2 & A_3 & A_4 & A_5 & A_6 & A_7 & B_2 & B_3 & B_4 & C_3 & C_4 & D_3 & D_4 & G_2 \\ \hline
e(\f s)&   2 &   3 &   4 &   5 &   6 &   7 &   8 &   4 &  6 &   8 &   6 &   8 &   4 &   7 &   6
\end{array}
\]
\end{example*}

\section{A Proof of a theorem in Penkov-Zuckerman}\label{sec_Penkov_Zuckerman}
This section is devoted to the proof of the following theorem:
\begin{thm}\label{thm_PWZ}  Let $\f k$ be an $\f{sl}_2$-subalgebra of a
semisimple Lie algebra $\f g$, none of whose simple factors is of
type $A1$. Then there exists a positive integer $b(\f k, \f g)$,
such that for every irreducible finite dimensional $\f g$-module
$V$, there exists an injection of $\f k$-modules $W\to V$, where
$W$ is an irreducible $\f k$-module of dimension less than $b(\f
k, \f g)$.
\end{thm}
\begin{proof}
We assume all of the structure of Sections \ref{subsec_Lie_theory}
and \ref{subsec_rep_theory} (ie: $\f h$, $\Phi$, $W$, etc.).
Consider a fixed $\f k$ and $\f g$.  For $k \in \bb N$, let $V(k)$
denote the irreducible, finite dimensional representation of $\f
k$ of dimension $k+1$. Each irreducible representation of $\f g$
may be regarded as a $\f k$ representation by restriction. As
before, for $\la \in P_+(\f g)$, we let $[L(\la), V(k)]$ denote
the multiplicity of $V(k)$ in $L(\la)$. Set:
\[g_0(\la) := \min \{ \dim V(k) |  k \in \bb N \text{ and } [L(\la),V(k)] \neq 0 \}.\]
For each fundamental weight $\omega_i$ ($1 \leq i \leq r$), let $m_i := m(\f g,\f k, i)$ (as in Definition \ref{dfn_m}).
By Proposition \ref{prop_m_positive}, $m_i \neq 0$ for all $i$.
Set $\delta_i := m_i \omega_i$, and define
$\mc C_0 := \{\sum_{i=1}^r a_i \omega_i| 0 \leq a_i < n_i  \}$.
We set:
\[b(\f k, \f g) := \max \left\{ g_0(\la) | \la \in \mc C_0\right\} + 1.\]
For $\q \in \bb N^r$, let $\mc C_\q := \left( \sum_{i=1}^r q_i
\delta_i\right) + \mc C_0$. By the division algorithm, the
collection of sets, $\left \{ \mc C_\q \; | \; \q  \in \bb N^r
\right \}$ partitions $P_+(\f g)$. We claim that for every $\q \in
\bb N^r$, $\max \{g_0(\la)|\la \in \mc C_q\} < b(\f k, \f g)$.

The result follows from this claim.  Indeed, let $V = L(\la)$, for
$\la \in P_+(\f g)$. There exists (a unique) $\q \in \bb N^r$ such
that $\la \in \mc C_\q$. Let $W$ be the irreducible $\f
k$-representation of dimension $g_0(\la)$.  By definition of
$g_0$, there exists an injection of $W$ into $L(\la)$.  By the
claim, $\dim W = g_0(\la) < b(\f k, \f g)$. The result follows.

We now will establish the claim by applying Proposition \ref{prop_semigroup}.

Let $\q \in \bb N^r$ and $\la^\prime \in \mc C_\q$. $\mu =
\sum_{i=1}^r q_i \delta_i$. Set: $\la := \la^\prime - \mu$. By
definition of $\mc C_\q$, we have $\la \in \mc C_0$, and so
$g_0(\la) < b(\f k, \f g)$. This means that there exists an
irreducible $\f k$-representation, $W$, such that $\dim W <
b(\f k, \f g)$ and $[W, L(\la)] \neq 0$.  By definition,
$\la^\prime = \la + \mu$. Clearly, \footnote{As a special case of
Proposition \ref{prop_semigroup}, we see that since
$L(\delta_i)^{\f k} \neq (0)$ for all $i$, we have $L(\mu)^{\f k}
\neq (0)$.} $L(\mu)^{\f k} \neq (0)$. Applying Proposition
\ref{prop_semigroup} we see that $[W, L(\la^\prime)] \neq 0$. This
means that $g_0(\la^\prime) \leq b(\f k, \f g)$.  The claim
follows.
\end{proof}

\begin{thm}
Let $\f k$ be a semisimple subalgebra of a
semisimple Lie algebra $\f g$, none of whose simple factors is in
the set $E(\f k)$.  Then there exists a positive integer $b(\f k, \f g)$,
such that for every irreducible finite dimensional $\f g$-module
$V$, there exists an injection of $\f k$-modules $W\to V$, where
$W$ is an irreducible $\f k$-module of dimension less than $b(\f
k, \f g)$.
\end{thm}
\begin{proof}
The proof is essentially the same as the proof of Theorem \ref{thm_PWZ}.
The only changes are a substitution of Proposition \ref{prop_m_positive_general}
for Proposition \ref{prop_m_positive}, and we index irreducible representation of
$\f k$ by $P_+(\f k)$ rather than $\bb N$.  We leave it to the reader to fill in the
details.
\end{proof}

The above theorems begs us to compute the smallest value of $b(\f
k, \f g)$. This is the subject of the Section \ref{sec_G2} for the
case when $\f g$ is the exceptional Lie algebra $G_2$ and $\f k$
is a principal $\f{sl}_2$-subalgebra of $\f g$.  Other examples
will follow in future work.

We remark that the number $b(\f k ,\f g)$ clearly depends on $\f
k$ (as the notation suggests).  Of course, there are only finitely
many $\f{sl}_2$-subalgebras in $\f g$, up to automorphism of $\f
g$.  We can therefore, consider the maximum value of $b(\f k, \f
g)$ as $\f k$ ranges over this finite set.  We will call this
number $b(\f g)$.  With this in mind, one might attempt to
estimate $b(\f g)$ for a given semisimple Lie algebra $\f g$.

On the other hand, there is a sense that one could fix $\f k$ to
be (say) a principal $\f sl_2$-subalgebra in some $\f g$.  We
could then consider the question of whether $b(\f k, \f g)$ is
bounded as $\f g$ varies (among semisimple Lie algebras with no
simple A1 factor).

Even more impressive would be allowing \emph{both} $\f g$ and $\f
k$ to vary.  It is certainly not clear that a bound would even
exists for $b(\f k, \f g)$.  If it did, we would be interested in
an estimate.

\section{Example: $G_2$}\label{sec_G2}

In this section we consider an example which illustrates the
result of Section \ref{sec_Penkov_Zuckerman}.

Let $G$ be a connected, simply connected, complex algebraic group
with $\f g \cong G_2$.  Let $K$ be a connected principal
$SL_2$-subgroup of $G$.  As before, we set $\f k = Lie(K)$.

We order the fundamental weights of $\f g$ so that $L(\omega_1) =
7$ and $\dim L(\omega_2) = 14$.  For the rest of this section, we
will refer to the representation $L(a \omega_1 + b \omega_2)$ (for
$a,b \in \bb N$) as $[a,b]$. In Table \ref{Table_dim_SL2_G2}, the
entry in row $i$ column $j$ is $\dim [i,j]^K$:
{\begin{table}[bht]$
\begin{array}{r|rrrrrrrrrrrrrrrrrrrr} i \backslash j
  & 0  & 1 & 2 &  3  & 4 & 5 & 6 & 7 & 8 & 9 & 10 & 11 & 12 & 13 & 14 & 15 & 16 & 17 & 18 & 19 \\ \hline
0 & 1 & 0 &  1  & 0  & 1  & 0  & 2  & 0  & 2  & 0  & 3  & 0  & 4  & 0  & 4  & 1  & 5  & 1   & 6  & 1    \\
1 & 0 & 0 &  0  & 1  & 0  & 1  & 1  & 1  & 2  & 2  & 2  & 3  & 3  & 4  & 4  & 5  & 5  & 6   & 7  & 7    \\
2 & 0 & 0 &  1  & 0  & 2  & 0  & 3  & 1  & 4  & 2  & 5  & 3  & 7  & 4  & 9  & 5  & 11 & 7   & 13 & 9    \\
3 & 0 & 0 &  0  & 1  & 1  & 2  & 2  & 3  & 3  & 5  & 5  & 6  & 7  & 8  & 9  & 11 & 11 & 13  & 14 & 16   \\
4 & 1 & 0 &  2  & 1  & 3  & 2  & 5  & 3  & 7  & 5  & 9  & 7  & 12 & 9  & 15 & 12 & 18 & 15  & 22 & 18   \\
5 & 0 & 1 &  1  & 2  & 2  & 4  & 4  & 6  & 6  & 8  & 9  & 11 & 12 & 14 & 15 & 18 & 19 & 22  & 23 & 26   \\
6 & 1 & 0 &  2  & 2  & 4  & 3  & 7  & 5  & 10 & 8  & 13 & 11 & 17 & 15 & 21 & 19 & 26 & 23  & 32 & 28   \\
7 &  0 & 1 &  1  & 3  & 3  & 5  & 6  & 8  & 9  & 12 & 12 & 16 & 17 & 20 & 22 & 25 & 27 & 31  & 33 & 37   \\
8 & 1 & 1 &  3  & 2  & 6  & 5  & 9  & 8  & 13 & 12 & 18 & 16 & 23 & 21 & 29 & 27 & 35 & 33  & 42 & 40   \\
9  & 0 & 1 &  2  & 4  & 4  & 7  & 8  & 11 & 12 & 16 & 17 & 21 & 23 & 27 & 29 & 34 & 36 & 41  & 44 & 49   \\
10 & 2 & 1 &  4  & 4  & 7  & 7  & 12 & 11 & 17 & 16 & 23 & 22 & 30 & 28 & 37 & 36 & 45 & 44  & 54 & 52   \\
11 & 0 & 2 &  2  & 5  & 6  & 9  & 10 & 14 & 16 & 20 & 22 & 27 & 29 & 35 & 37 & 43 & 46 & 52  & 56 & 62   \\
12 & 2 & 2 &  5  & 5  & 9  & 9  & 15 & 14 & 21 & 21 & 28 & 28 & 37 & 36 & 46 & 45 & 56 & 55  & 67 & 66   \\
13 & 0 & 2 &  3  & 6  & 7  & 11 & 13 & 17 & 19 & 25 & 27 & 33 & 36 & 42 & 46 & 53 & 56 & 64  & 68 & 76   \\
14 & 2 & 2 &  6  & 6  & 11 & 11 & 17 & 18 & 25 & 25 & 34 & 34 & 44 & 44 & 55 & 55 & 67 & 67  & 80 & 80   \\
15 & 1 & 3 &  4  & 8  & 9  & 14 & 16 & 21 & 24 & 30 & 33 & 40 & 44 & 51 & 55 & 64 & 68 & 77  & 82 & 91   \\
16 & 3 & 3 &  7  & 8  & 13 & 14 & 21 & 21 & 30 & 31 & 40 & 41 & 52 & 53 & 65 & 66 & 79 & 80  & 95 & 95   \\
17 & 0 & 4 &  5  & 9  & 11 & 16 & 19 & 25 & 28 & 35 & 39 & 47 & 51 & 60 & 65 & 74 & 80 & 90  & 96 & 107  \\
18 & 3 & 3 &  8  & 9  & 15 & 16 & 24 & 25 & 34 & 36 & 46 & 48 & 60 & 61 & 75 & 77 & 91 & 93  & 109& 111  \\
19 & 1 & 4 &  6  & 11 & 13 & 19 & 22 & 29 & 33 & 41 & 45 & 54 & 60 & 69 & 75 & 86 & 92 & 104 & 111& 123  \\
\end{array}$
\caption{$\dim [i,j]^K$ where $K$ is a principle $SL_2$-subgroup in $G_2$.}\label{Table_dim_SL2_G2}
\end{table}

%\begin{center}
%\includegraphics[75,0][250,500]{G2_Picture.pdf}\\
%The above picture is a $40 \times 40$ version of Table
%\ref{Table_dim_SL2_G2} with color coded multiplicity.
%\end{center}

\medskip Table \ref{Table_dim_SL2_G2} was generated by an
implementation of the Weyl character formula (see Theorem
\ref{thm_Weyl_character_formula}) for the group $G_2$ using the
computer algebra system MAPLE$^\circledR$. The characters were
restricted to a maximal torus in $K$, thus allowing us to find the
character of $[i,j]$ as a representation of $SL_2(\bbC)$.  This
character was used to compute the dimension of the invariants for
$K$.

Using the same implementation we can compute the values of
$g_0(\la)$ for the pair $(\f k, \f g)$.  We display these data in
the table below for $0 \leq i,j \leq 19$.
\begin{table}[bht]$
\begin{array}{r|rrrrrrrrrrrrrrrrrrrr}
i \backslash j
   &      0  &  1 &   2 &   3 &   4 &   5 &   6 &   7 &   8 &   9 &  10 &  11 &  12 &  13 &  14  & 15 &  16 &  17 &  18 &  19   \\ \hline
0  &      1  &  3 &   1 &   3 &   1 &   3 &   1 &   3 &   1 &   3 &   1 &   3 &   1 &   3 &   1  &  1 &   1 &   1 &   1 &   1   \\
1  &      7  &  5 &   3 &   1 &   3 &   1 &   1 &   1 &   1 &   1 &   1 &   1 &   1 &   1 &   1  &  1 &   1 &   1 &   1 &   1   \\
2  &      5  &  3 &   1 &   3 &   1 &   3 &   1 &   1 &   1 &   1 &   1 &   1 &   1 &   1 &   1  &  1 &   1 &   1 &   1 &   1   \\
3  &      3  &  3 &   3 &   1 &   1 &   1 &   1 &   1 &   1 &   1 &   1 &   1 &   1 &   1 &   1  &  1 &   1 &   1 &   1 &   1   \\
4  &      1  &  3 &   1 &   1 &   1 &   1 &   1 &   1 &   1 &   1 &   1 &   1 &   1 &   1 &   1  &  1 &   1 &   1 &   1 &   1   \\
5  &      3  &  1 &   1 &   1 &   1 &   1 &   1 &   1 &   1 &   1 &   1 &   1 &   1 &   1 &   1  &  1 &   1 &   1 &   1 &   1   \\
6  &      1  &  3 &   1 &   1 &   1 &   1 &   1 &   1 &   1 &   1 &   1 &   1 &   1 &   1 &   1  &  1 &   1 &   1 &   1 &   1   \\
7  &      3  &  1 &   1 &   1 &   1 &   1 &   1 &   1 &   1 &   1 &   1 &   1 &   1 &   1 &   1  &  1 &   1 &   1 &   1 &   1   \\
8  &      1  &  1 &   1 &   1 &   1 &   1 &   1 &   1 &   1 &   1 &   1 &   1 &   1 &   1 &   1  &  1 &   1 &   1 &   1 &   1   \\
9  &      3  &  1 &   1 &   1 &   1 &   1 &   1 &   1 &   1 &   1 &   1 &   1 &   1 &   1 &   1  &  1 &   1 &   1 &   1 &   1   \\
10 &      1  &  1 &   1 &   1 &   1 &   1 &   1 &   1 &   1 &   1 &   1 &   1 &   1 &   1 &   1  &  1 &   1 &   1 &   1 &   1   \\
11 &      3  &  1 &   1 &   1 &   1 &   1 &   1 &   1 &   1 &   1 &   1 &   1 &   1 &   1 &   1  &  1 &   1 &   1 &   1 &   1   \\
12 &      1  &  1 &   1 &   1 &   1 &   1 &   1 &   1 &   1 &   1 &   1 &   1 &   1 &   1 &   1  &  1 &   1 &   1 &   1 &   1   \\
13 &      3  &  1 &   1 &   1 &   1 &   1 &   1 &   1 &   1 &   1 &   1 &   1 &   1 &   1 &   1  &  1 &   1 &   1 &   1 &   1   \\
14 &      1  &  1 &   1 &   1 &   1 &   1 &   1 &   1 &   1 &   1 &   1 &   1 &   1 &   1 &   1  &  1 &   1 &   1 &   1 &   1   \\
15 &      1  &  1 &   1 &   1 &   1 &   1 &   1 &   1 &   1 &   1 &   1 &   1 &   1 &   1 &   1  &  1 &   1 &   1 &   1 &   1   \\
16 &      1  &  1 &   1 &   1 &   1 &   1 &   1 &   1 &   1 &   1 &   1 &   1 &   1 &   1 &   1  &  1 &   1 &   1 &   1 &   1   \\
17 &      3  &  1 &   1 &   1 &   1 &   1 &   1 &   1 &   1 &   1 &   1 &   1 &   1 &   1 &   1  &  1 &   1 &   1 &   1 &   1   \\
18 &      1  &  1 &   1 &   1 &   1 &   1 &   1 &   1 &   1 &   1 &   1 &   1 &   1 &   1 &   1  &  1 &   1 &   1 &   1 &   1   \\
19 &      1  &  1 &   1 &   1 &   1 &   1 &   1 &   1 &   1 &   1
&   1 &   1 &   1 &   1 &   1  &  1 &   1 &   1 &   1 &   1
\end{array}$
\caption{Value of $g_0(i \omega_1 + j \omega_2)$ for a principal $SL_2$-subgroup in $G_2$.}\label{Table_g0_SL2G2}
\end{table}

Of particular interest is the 7 in row 1 column 0. The irreducible
$G_2$-representation $[1,0]$ is irreducible when restricted to a
principal $\f{sl}_2$ subalgebra. And therefore,
$g_0(\omega_1)=\dim [1,0] = 7$.  Note that most entries in Table
\ref{Table_g0_SL2G2} are 1.  Following the proof in Section
\ref{sec_Penkov_Zuckerman} we see that $b(\f k, \f g) = 7 + 1 =
8$. That is, every finite dimensional representation of $G_2$
contains an irreducible $\f k$-representation of dimension less
than 8.

Even more interesting is the fact that Table \ref{Table_g0_SL2G2}
suggests that there are only 26 ordered pairs ($a,b$) such that
$g_0(a \omega_1 + b \omega_2) > 1$.  This is indeed the case:
\begin{thm}\label{thm_G2_invariants}
For all $a,b \in \bb N$, $[a,b]^K \neq (0)$ except for
the following list of 26 exceptions:
\begin{eqnarray*}
&[0, 1], [0, 3], [0, 5], [0, 7], [0, 9], [0, 11], [ 0,13], [1, 0], [1, 1], \\
&[1, 2], [1, 4], [2, 0], [2, 1], [2, 3], [2,  5], [ 3, 0], [3, 1], [3, 2], \\
&[4, 1], [5, 0], [6, 1], [7, 0], [9, 0], [11, 0], [13, 0], [17,0].
\end{eqnarray*}
\end{thm}
\begin{proof}
By inspection of Table \ref{Table_dim_SL2_G2} (or Table
\ref{Table_g0_SL2G2}), it suffices to show that for all $a,b \in
\bb N$ such that $a >24$ or $b>24$ we have $\dim [a,b]^K > 0$. Let
$\mu_1, \cdots, \mu_6$ denote the highest weights of the
representations $[0,2],[0,17],[4,0], [15,0],[5,1]$, and  $[1,3]$.
From the table we see that each of these representations has a
$K$-invariant. Let $Z = \sum_{i=1}^6 \bb N \mu_i \subset P_+(\f
g)$.  By Proposition \ref{prop_semigroup}, (for $W$ trivial), each
element of $Z$ has a $K$-invariant.  It is easy to see that $E =
P_+(\f g) - Z$ is a finite set.  In fact, computer calculations show
that $E$ consists of the following 194 elements:

\noindent{\tiny [0, 1], [0, 3], [0, 5], [0, 7], [0, 9], [0, 11],
[0, 13], [0, 15], [1, 0], [1, 1], [1, 2], [1, 4], [1, 6], [1, 8],
[1, 10], [1, 12], [1, 14], [1, 16], [1, 18], [2, 0], [2, 1], [2,
2], [2, 3], [2, 4], [[0, 1], [0, 3], [0, 5], [0, 7], [0, 9], [0,
11], [0, 13], [0, 15], [1, 0], [1, 1], [1, 2], [1, 4], [2, 0], [2,
1], [2, 2], [2, 3], [2, 4], [2, 5], [2, 7], [3, 0], [3, 1], [3,
2], [3, 3], [3, 4], [3, 5], [3, 6], [3, 7], [3, 8], [3, 10], [4,
1], [4, 3], [4, 5], [4, 7], [4, 9], [4, 11], [4, 13], [5, 0], [5,
2], [5, 4], [6, 1], [6, 3], [6, 5], [7, 0], [7, 2], [7, 4], [8,
1], [8, 3], [8, 5], [9, 0], [9, 2], [9, 4], [10, 1], [10, 3], [10,
5], [11, 0], [11, 2], [11, 4], [12, 1], [12, 3], [12, 5], [13, 0],
[13, 2], [13, 4], [14, 1], [14, 3], [14, 5], [16, 1], [17, 0],
[17, 2], [17, 4], [18, 1], [18, 3], [18, 5], [2, 5], [2, 7], [2,
9], [2, 11], [2, 13], [2, 15], [2, 17], [2, 19], [2, 21], [3, 0],
[3, 1], [3, 2], [3, 3], [3, 4], [3, 5], [3, 6], [3, 7], [3, 8],
[3, 10], [3, 12], [3, 14], [3, 16], [3, 18], [3, 20], [3, 22], [3,
24], [4, 1], [4, 3], [4, 5], [4, 7], [4, 9], [4, 11], [4, 13], [4,
15], [5, 0], [5, 2], [5, 4], [5, 6], [5, 8], [5, 10], [5, 12], [5,
14], [5, 16], [6, 0], [6, 1], [6, 2], [6, 3], [6, 5], [6, 7], [6,
9], [6, 11], [6, 13], [6, 15], [6, 17], [6, 19], [7, 0], [7, 1],
[7, 2], [7, 3], [7, 4], [7, 5], [7, 6], [7, 8], [7, 10], [7, 12],
[7, 14], [7, 16], [7, 18], [7, 20], [7, 22], [8, 1], [8, 3], [8,
5], [8, 7], [8, 9], [8, 11], [8, 13], [8, 15], [9, 0], [9, 2], [9,
4], [9, 6], [9, 8], [9, 10], [9, 12], [9, 14], [9, 16], [10, 0],
[10, 1], [10, 3], [10, 5], [10, 7], [10, 9], [10, 11], [10, 13],
[10, 15], [10, 17], [11, 0], [11, 1], [11, 2], [11, 3], [11, 4],
[11, 6], [11, 8], [11, 10], [11, 12], [11, 14], [11, 16], [11,
18], [11, 20], [12, 1], [12, 3], [12, 5], [12, 7], [12, 9], [12,
11], [12, 13], [12, 15], [13, 0], [13, 2], [13, 4], [13, 6], [13,
8], [13, 10], [13, 12], [13, 14], [13, 16], [14, 0], [14, 1], [14,
3], [14, 5], [14, 7], [14, 9], [14, 11], [14, 13], [14, 15], [14,
17], [15, 1], [16, 1], [17, 0], [17, 2], [17, 4], [18, 0], [18,
1], [18, 3], [18, 5], [18, 7], [19, 1], [21, 0], [21, 2], [22, 0],
[22, 1], [22, 3], [22, 5], [23, 1], [25, 0], [26, 0], [26, 1],
[26, 3], [27, 1], [29, 0], [30, 1], [31, 1], [33, 0], [34, 1],
[37, 0], [38, 1], [41, 0], [42, 1], [46, 1]}.

If one checks, we see that each of these has an invariant except for the 26 exceptional values in the statement of the theorem.

Unfortunately, our tables (in this paper) are not big enough to
see all of these weights.  For this reason, we can alternately
consider the set $Z^\prime = \sum_{i=1}^9 \bb N \mu^\prime_i$,
where $\mu_1^\prime, \cdots, \mu_9^\prime$ are the highest weights
of the representations $[0,2], [0,17], [4,0], [6, 0],
[15,0],[5,1], [1,3],[7,1],$ and $[1,6]$. Again, each of these has
a $K$-invariant (from the tables) and therefore by Proposition
\ref{prop_semigroup}, each element of $Z^\prime$ has a
$K$-invariant.  Set $E^\prime = P_+(\f g) - Z^\prime$. $E^\prime$
has 73 elements.  They are: {\tiny [0, 1], [0, 3], [0, 5], [0, 7],
[0, 9], [0, 11], [0, 13], [0, 15], [1, 0], [1, 1], [1, 2], [1, 4],
[2, 0], [2, 1], [2, 2], [2, 3], [2, 4], [2, 5], [2, 7], [3, 0],
[3, 1], [3, 2], [3, 3], [3, 4], [3, 5], [3, 6], [3, 7], [3, 8],
[3, 10], [4, 1], [4, 3], [4, 5], [4, 7], [4, 9], [4, 11], [4, 13],
[5, 0], [5, 2], [5, 4], [6, 1], [6, 3], [6, 5], [7, 0], [7, 2],
[7, 4], [8, 1], [8, 3], [8, 5], [9, 0], [9, 2], [9, 4], [10, 1],
[10, 3], [10, 5], [11, 0], [11, 2], [11, 4], [12, 1], [12, 3],
[12, 5], [13, 0], [13, 2], [13, 4], [14, 1], [14, 3], [14, 5],
[16, 1], [17, 0], [17, 2], [17, 4], [18, 1], [18, 3], [18, 5]}.\\
This time, each of these representations is on our tables. And so,
the reader may (and should) check that, apart from the 26
exceptional values in the statement of the theorem, all of these
have a positive dimension of $K$-invariants.
\end{proof}
\section{Tables for maximal parabolic subalgebras.}\label{sec_tables} $\empty$
In the proof of Corollary \ref{cor_SL2_invariants_noA2} we
needed to see that if $\f g$ was a simple Lie algebra not of type
A2, then $\dim (\f g/ [\f p, \f p]) > 3$ for a maximal parabolic
subalgebra $\f p$.  If $\f p = \f l \oplus \f u^+$, $[\f l, \f l]
= \f l_{ss}$ then we can deduce that this inequality is equivalent
to $\dim (\f g / \f l_{ss}) > 5$ (using the formula $\dim (\f g/ [\f
p, \f p]) = \frac{1}{2}\left( \dim \f g/ \f l_{ss} + 1\right)$. It
is possible to exhaust over all possible cases to establish this
fact.  The exceptional cases are in the following table:

\begin{center}
Dimension of $\f g / \f l_{ss}$ for the exceptional groups.\\
{\tiny $\left.\begin{tabular}{|c||r|r|r|r|r|r|r|r|} \hline
  G $\setminus$ $\f l_{ss}$ & 1 & 2 & 3 & 4 & 5 & 6 & 7 & 8 \\
  \hline \hline
  G2     & A1 & A1 \\ \cline{1-3}
$\dim \f g / \f l_{ss}$
         & 11 & 11 \\ \cline{1-5} \cline{1-5}
  F4     & C3 & A2A1 & A2A1 & B3  \\ \cline{1-5}
         & 31 & 41 & 41 & 31  \\ \cline{1-7}
  E6     & D5 & A5 & A4A1 & A2A2A1 & A4A1 & D5 \\ \cline{1-7}
         & 33 & 43 & 51 & 59 & 51 & 33 \\ \cline{1-8}
  E7     & D6 & A6 & A1A5 & A1A2A3 & A2A4 & D5A1 & E6 \\  \cline{1-8}
         & 67 & 85 & 95 & 107 & 101 & 85 & 55 \\  \hline
  E8     & D7 & A7 & A1A6 & A1A2A4 & A4A3 & D5A2 & E6A1 & E7 \\ \hline
         & 157 & 185 & 197 & 213 & 209 & 195 & 167 & 115 \\ \hline
\end{tabular} \right |$}\end{center}

For the classical cases, the situation reduces to the examination
of several families of parabolic subalgebras.  In each family, we
can determine the dimension of $\f g/ \f l_{ss}$ for any of
parabolic $\f p = \f l \oplus \f u^+$, from the formulas:
\begin{eqnarray*}
\dim A_n &=& (n+1)^2-1 \\
\dim B_n &=&  n(2n+1)  \\
\dim C_n &=&  n(2n+1)  \\
\dim D_n &=&  n(2n-1)
\end{eqnarray*}
For example, if $\f g=A_n$ $(n \geq 2)$ and $\f l_{ss}=A_p \oplus
A_{q}$ for $p+q=n-1$ (set $\dim A_0 := 0$) we have $\dim (\f g/ \f
l_{ss}) = \dim A_n - \dim A_p - \dim A_{q} =
n^2-p^2-q^2+2n+2p+2q-2$.  Upon examination of the possible values
of this polynomial over the parameter space we see that the
smallest value is 5.  This case occurs only for $n=2$.  All other
values of $p$ and $q$ give rise to larger values. In fact, for $\f
g$ classical and not of type A the dimension is also greater than
5. The low rank cases are summarized in the following tables:

Dimension of $\f g / \f l_{ss}$ for type A. \hspace{3cm} Dimension of $\f g / \f l_{ss}$ for type D. \\
{\tiny $\left.\begin{tabular}{|c||r|r|r|r|r|} \hline
  G $\setminus$ $\f l_{ss}$ & 1 & 2 & 3 & 4 & 5 \\ \hline \hline
A2 &   A1 &   A1  \\ \cline{1-3}
   &    5 &    5  \\ \cline{1-4}
A3 &   A2 & A1A1 & A2 \\ \cline{1-4}
   &    7 &    9 &  7 \\ \cline{1-5}
A4 &   A3 & A1A2 & A1A2 & A3 \\ \cline{1-5}
   &    9 &   13 &   13 &  9 \\ \cline{1-6}
A5 &   A4 & A1A3 & A2A2 & A1A3 & A4 \\ \cline{1-6}
   &   11 &   17 &   19 &   17 & 11 \\ \cline{1-6} \hline
\end{tabular}\right|$ \hspace{1cm}
$\left.
\begin{tabular}{|c||r|r|r|r|r|r|} \hline
  G $\setminus$ $\f l_{ss}$ & 1 & 2 & 3 & 4 & 5 & 6 \\ \hline \hline
D3 & A1A1 &   A2 &   A2 \\ \cline{1-4}
   &    9 &    7 &    7 \\ \cline{1-5}
D4 &   A3 &A1A1A1&   A3 &   A3 \\ \cline{1-5}
   &   13 &   19 &   13 &   13 \\ \cline{1-6}
D5 &   D4 & A1D3 &A2A1A1&   A4 & A4 \\ \cline{1-6}
   &   17 &   27 &   31 &   21 & 21 \\ \cline{1-6} \hline
D6 &   D5 & A1D4 & A2D3 &A1A1A3& A5 & A5 \\ \cline{1-7}
   &   21 &   35 &   43 &   45 & 31 & 31 \\ \cline{1-7} \hline
\end{tabular} \right |$}

\medskip \noindent
Dimension of $\f g / \f l_{ss}$ for type B.  \hspace{3cm} Dimension of $\f g / \f l_{ss}$ for type C.\\
{\tiny $\left.\begin{tabular}{|c||r|r|r|r|r|} \hline
  Group $\setminus$ $\f l_{ss}$ & 1 & 2 & 3 & 4 & 5 \\ \hline \hline
B2 &   A1 &   A1  \\ \cline{1-3}
   &    7 &    7  \\ \cline{1-4}
B3 &   B2 & A1A1 &   A2 \\ \cline{1-4}
   &   11 &   15 &   13 \\ \cline{1-5}
B4 &   B3 & A1B2 & A1A2 &   A3 \\ \cline{1-5}
   &   15 &   23 &   25 &   21 \\ \cline{1-6}
B5 &   B4 & A1B3 & A2B2 & A1A3 & A4 \\ \cline{1-6}
   &   19 &   31 &   37 &   37 & 31 \\ \cline{1-6} \hline
\end{tabular} \right |$ \hspace{10mm}
$\left.
\begin{tabular}{|c||r|r|r|r|r|} \hline
  Group $\setminus$ $\f l_{ss}$ & 1 & 2 & 3 & 4 & 5 \\ \hline \hline
C2 &   A1 &   A1  \\ \cline{1-3}
   &    7 &    7  \\ \cline{1-4}
C3 &   C2 & A1A1 &   A2 \\ \cline{1-4}
   &   11 &   15 &   13 \\ \cline{1-5}
C4 &   C3 & A1C2 & A1A2 &   A3 \\ \cline{1-5}
   &   15 &   23 &   25 &   21 \\ \cline{1-6}
C5 &   C4 & A1C3 & A2C2 & A1A3 & A4 \\ \cline{1-6}
   &   19 &   31 &   37 &   37 & 31 \\ \cline{1-6} \hline
\end{tabular} \right |$}

%\clearpage
%\begin{figure}[htbp]
%        \includegraphics[10,20][350,600]{Picture_G2.JPG}
%\caption{Table \ref{Table_dim_SL2_G2} with a $40 \times 40$ grid
%and colored multiplicity.} \label{Picture_G2}
%\end{figure}
%\clearpage

%%\bibliographystyle{plain}
%\bibliographystyle{alpha}
%\bibliography{citations}

\begin{thebibliography}{Dyn52b}

\bibitem[Bor91]{Borel-LAG}
Armand Borel.
\newblock {\em Linear algebraic groups}, volume 126 of {\em Graduate Texts in
  Mathematics}.
\newblock Springer-Verlag, New York, second edition, 1991.

\bibitem[Bou89]{Bourbaki-Lie123}
Nicolas Bourbaki.
\newblock {\em Lie groups and {L}ie algebras. {C}hapters 1--3}.
\newblock Elements of Mathematics (Berlin). Springer-Verlag, Berlin, 1989.
\newblock Translated from the French, Reprint of the 1975 edition.

\bibitem[Bou02]{Bourbaki-Lie456}
Nicolas Bourbaki.
\newblock {\em Lie groups and {L}ie algebras. {C}hapters 4--6}.
\newblock Elements of Mathematics (Berlin). Springer-Verlag, Berlin, 2002.
\newblock Translated from the 1968 French original by Andrew Pressley.

\bibitem[CM93]{CM}
David~H. Collingwood and William~M. McGovern.
\newblock {\em Nilpotent orbits in semisimple {L}ie algebras}.
\newblock Van Nostrand Reinhold Mathematics Series. Van Nostrand Reinhold Co.,
  New York, 1993.

\bibitem[Dyn52a]{Dynkin-maximal}
E.~B. Dynkin.
\newblock Maximal subgroups of the classical groups.
\newblock {\em Trudy Moskov. Mat. Ob\v s\v c.}, 1:39--166, 1952.

\bibitem[Dyn52b]{Dynkin-semisimple}
E.~B. Dynkin.
\newblock Semisimple subalgebras of semisimple {L}ie algebras.
\newblock {\em Mat. Sbornik N.S.}, 30(72):349--462 (3 plates), 1952.

\bibitem[Eis95]{Eisenbud}
David Eisenbud.
\newblock {\em Commutative algebra, with a view toward algebraic geometry},
  volume 150 of {\em Graduate Texts in Mathematics}.
\newblock Springer-Verlag, New York, 1995.

\bibitem[FH91]{FH}
William Fulton and Joe Harris.
\newblock {\em Representation theory}.
\newblock Springer-Verlag, New York, 1991.
\newblock A first course, Readings in Mathematics.

\bibitem[GW98]{GW}
R.~Goodman and N.R. Wallach.
\newblock {\em Representations and invariants of the classical groups}.
\newblock Cambridge University Press, Cambridge, 1998.

\bibitem[Har77]{Hartshorne}
Robin Hartshorne.
\newblock {\em Algebraic geometry}.
\newblock Springer-Verlag, New York, 1977.
\newblock Graduate Texts in Mathematics, No. 52.

\bibitem[Kna02]{Knapp-Beyond}
Anthony~W. Knapp.
\newblock {\em Lie groups beyond an introduction}, volume 140 of {\em Progress
  in Mathematics}.
\newblock Birkh\"auser Boston Inc., Boston, MA, second edition, 2002.

\bibitem[Kos59]{Kostant}
Bertram Kostant.
\newblock The principal three-dimensional subgroup and the {B}etti numbers of a
  complex simple {L}ie group.
\newblock {\em Amer. J. Math.}, 81:973--1032, 1959.

\bibitem[Pop74]{Popov-UFD}
V.~L. Popov.
\newblock Picard groups of homogeneous spaces of linear algebraic groups and
  one-dimensional homogeneous vector fiberings.
\newblock {\em Izv. Akad. Nauk SSSR Ser. Mat.}, 38:294--322, 1974.

\bibitem[PV89]{PV}
V.~L. Popov and {\`E}.~B. Vinberg.
\newblock Invariant theory.
\newblock In {\em Algebraic geometry, 4 (Russian)}, Itogi Nauki i Tekhniki,
  pages 137--314, 315. Akad. Nauk SSSR Vsesoyuz. Inst. Nauchn. i Tekhn.
  Inform., Moscow, 1989.

\bibitem[PZ04]{PZ}
Ivan Penkov and Gregg Zuckerman.
\newblock {Generalized Harish-Chandra Modules: A New Direction}.
\newblock {\em Acta Applicandae Mathematicae}, 81(1):311--326, 2004.

\bibitem[Sha94]{AGIV}
I.~R. Shafarevich, editor.
\newblock {\em Algebraic geometry. {IV}}, volume~55 of {\em Encyclopaedia of
  Mathematical Sciences}.
\newblock Springer-Verlag, Berlin, 1994.
\newblock Linear algebraic groups. Invariant theory, A translation of {\it
  Algebraic geometry. }4 (Russian), Akad.\ Nauk SSSR Vsesoyuz.\ Inst.\ Nauchn.\
  i Tekhn.\ Inform., Moscow, 1989 [MR 91k:14001], Translation edited by A. N.
  Parshin and I. R. Shafarevich.

\bibitem[Spr89]{Springer}
T.~A. Springer.
\newblock Linear algebraic groups.
\newblock In {\em Algebraic geometry, 4 (Russian)}, Itogi Nauki i Tekhniki,
  pages 5--136, 310--314, 315. Akad. Nauk SSSR Vsesoyuz. Inst. Nauchn. i Tekhn.
  Inform., Moscow, 1989.
\newblock Translated from the English.

\bibitem[Vos69]{Vosk-UFD}
V.~E. Voskresenski\u{\i}.
\newblock Picard groups of linear algebraic groups.
\newblock In {\em Studies in Number Theory, No. 3 (Russian)}, pages 7--16.
  Izdat. Saratov. Univ., Saratov, 1969.

\bibitem[Vos98]{Vosk}
V.~E. Voskresenski\u{\i}.
\newblock {\em Algebraic groups and their birational invariants}, volume 179 of
  {\em Translations of Mathematical Monographs}.
\newblock American Mathematical Society, Providence, RI, 1998.
\newblock Translated from the Russian manuscript by Boris Kunyavski [Boris
  \`E.\ Kunyavski\u\i].

\bibitem[VP72]{VP1972}
{\`E}.~B. Vinberg and V.~L. Popov.
\newblock A certain class of quasihomogeneous affine varieties.
\newblock {\em Izv. Akad. Nauk SSSR Ser. Mat.}, 36:749--764, 1972.

\end{thebibliography}

\def\cprime{$'$} \def\cprime{$'$}

\end{document}